\documentclass[]{article}

\newtheorem{lemma}{Lemma}[section]
\newtheorem{cor}{Corollary}[section]
\newtheorem{prop}{Proposition}[section]

\newtheorem{thm}{Theorem}[section]

\newcommand{\Co}{ \mathbf{C}}
\newcommand{\R}{ {\mathbf R}}
\newcommand{\U}{{\mathcal U}}
\newcommand{\pf}{\textsc {Proof }}

\newcommand{\Diff}{  \textrm {Diff} }
\newcommand{\Tr}{  \textrm  {Tr}}

\newcommand\ip[2]{   \langle#1 , #2 \rangle} 

\begin{document}

\title{ Factor Representations of Diffeomorphism Groups}
\author{Robert Boyer}
\maketitle

\begin{abstract}
General semifinite factor representations
of the diffeomorphism group of euclidean space
are constructed by means of a canonical correspondence with the finite factor representations
of the inductive limit unitary group.  
This construction includes  the
quasi-free representations of the canonical commutation and
anti-commutation relations.
To establish this correspondence requires a non-linear form
of complete positivity as developed by Arveson.
We also compare the asymptotic character formula for the
unitary group with the thermodynamic ($N/V$) limit
construction for diffeomorphism group representations.
\end{abstract}

\section{Introduction}

The purpose of this work is to explore  connections between the representation theory 
of inductive limit classical groups and the representation theory of the group $\Diff_c (\R^d)$
of compactly supported diffeomorphisms.

We need to extend the method of   ``holomorphic extension'' of Ol'shanskii \cite{O2}
from the full unitary group to certain inductive limit groups.  Our approach also uses
a non-linear form of complete positivity given by Arveson \cite{A}.  
The corresponding representations of the diffeomorphism group are semifinite
factor representations, typically, of type $\textrm{II}_\infty$.
Further,  certain subfamilies of these representations may be described in terms of
the quasi-free representations of the commutation relations
\cite{SV}. The overall framework for constructing  representations of the diffeomorphism 
group as inductive limits was first described in \cite{MS}.

See \cite{G3} for a survey of the representation theory of the diffeomorphism group.

\section{Diffeomorphism Group and Its Standard Representation}

Let $X$ denote the euclidean space $\R^d$,
where $d \geq 2$.
Let $Y$ denote an open subset of
$X$.
We let $\Diff_c(Y)$ denote the group of all compactly
supported diffeomorphisms of $Y$ with the usual
topology.
That is, a diffeomorphism $\psi$ has  \textit {compact support}
if it is equal to the identity on the complement of
some relatively compact open subset of $X$.
Further, a sequence of compactly supported
diffeomorphisms $\{ \psi_n \}$ converges to
$\psi$ if there is a common relatively compact
open set on whose complement all diffeomorphisms
are equal to the identity and, on the open set itself,
the diffeomorphisms converge in the ${\mathcal {C}}^\infty$
topology \cite{VGG}.

Let ${\cal F}$ denote the increasing sequence
$X_1 \subset X_2 \subset \dots$ 
 of connected open subsets of $X$,
such that $X_N$ is a proper subset of $X_{N+1}$,
where we allow the possibility that $m(X_1)$ is
infinite.
Then $\Diff_c(X) = \bigcup \Diff_c(X_N)$.

The usual orthogonal action of
$\Diff_c(X_N)$ on the real Hilbert space $L^2_{\R}(X_N),$ 
formed relative to Lebesgue measure $m,$
gives rise to the standard representation $T$ of the 
diffeomorphism group on the complex Hilbert space
$L^2(X_N)$, where 

\[
[T(\psi)f](x) = J_\psi(x)^{1/2}(x) f( \psi^{-1}x),
\quad \textrm{where  }  f  \in L^2(X_N), 
\]

\noindent
and $J_{\psi}(x) $ denotes the Radon-Nikodym derivative
$ dm( \psi^{-1}x)/dm$.

\section{Inductive Limit Groups and Their Representations}
\label{sec:limitgroup}

Let $L$ be a real separable infinite-dimensional Hilbert space with
an increasing sequence ${\cal L}$ of closed subspaces $\{L_N\}_{N=0}^\infty$.
We let ${\cal  {O}} = {\cal {O}}_{\cal {L}} = \lim_{\to}  O(L_N)$ 
be the inductive limit group where each
$O(L_N)$ is the usual orthogonal group with the strong operator topology.
Then the embedding $T$ given above gives a natural map from $\Diff_c(X)$
into ${\cal  {O}}_{\cal L}$ where $L_N = L^2(X_N)$ and $\Diff_c(X_N)$ maps into
$O(L_N)$.

We make the convention that by a representation of a topological group we mean a 
strongly continuous unitary representation.

Let $H_N$ be the complexification of $L_N$. We form $\U = \lim_{\to} U(H_N)$
so ${\cal  O}$ is a natural subgroup of $\U$.  Recall that a strongly continuous 
unitary representation of $U(H_N)$ is called {\it holomorphic} if its decomposition
into irreducibles whose signatures are all non-negative.  Now if $\pi_N$ is a representation
of $U(H_N)$ then its restriction to $O(L_N)$ has the same commutant. In particular,
restriction preserves irreducibility \cite{O2}, p 24.  Then by elementary reasoning
about generating nests of von Neumann algebras,
we find that
if $\pi$ is a representation of $\U$, then its restriction to ${\cal {O}}$ has the same commutant.
Moreover,  by the results of \cite{VGG}, the further restriction of $\pi$ to the image under
$T$ of $\Diff_c(X)$ has the same commutant.  In particular, the restriction of a factor
representation  $\pi$of $\U$  to $T[\Diff_c(X)]$ is factor.

We now indicate how  the KMS positive definite functions introduced in \cite{SV}
extend continuously to $\U$ and are semifinite factor representations.  In general, it
is a difficult problem to decide if a such a KMS positive definite function is factor, see 
\cite{SV, B1}.

The holomorphic finite characters of $U(\infty)$ are naturally parametrized by
 $\cal M$,  a class of mereomorphic functions,  where the character
$\chi_f(V) = \det [ f(V) ]$ is given by
$f \in {\cal {M}}$ with $f(z)$ in the form:
$e^{\lambda(z-1)} 
\prod_{k=1}^\infty (1-\beta_k+\beta_k z) / (1+\alpha_k - \alpha_k z)$,
where $0 \leq \lambda, \alpha_k, \beta_k$ and 
$\sum_{k=1}^\infty (\alpha_k  + \beta_k) < \infty$ \cite{B1,VK}.

Given the increasing sequence $\{H_N\}_{N=1}^\infty$, we choose unit vectors
$f_N \in H_N \ominus H_{N-1}, \, N \geq 1$. 
Note that the two collections of vectors are mutually orthornormal.
Let $H_-$ be the closed Hilbert space generated by $\langle f_1, f_2, \dots \rangle$
and $H_+ = H_-^\perp$.

As a convention, we now fix $U(\infty)$ to be the inductive limit unitary group
formed relative to the subspaces generated by 
$\langle f_1, \dots, f_N \rangle$ and $U(2 \infty)$ relative to
$\langle e_1, \dots, e_N, f_1, \dots, f_N \rangle$. Let $F$ denote the natural
projection of $H$ onto $H_-$. 

The function $\phi_{F,f}(V) = \det[ f(FVF)]$, for $V \in U(2\infty)$, is positive definite.
It extends continuously to $\U$ and then through restriction gives a positive definite
function on $\Diff_c(X)$. Let $\pi_{F,f}$ denote the corresponding cyclic representation.

\begin{thm} \label{thm:main}
$\pi_{F,f}$ is a semifinite factor representation of  $\Diff_c(X)$.
It  is irreducible if and only if
$f(z)=z^m$,
otherwise it is a $\textrm{II}_\infty$ 
factor representation.
\end{thm}

By \cite{VGG}, theorem 1.2, if $m(X_1) = + \infty$, the result of our theorem
holds for $\textrm{SDiff}_c(X)$, the subgroup of measure-preserving
diffeomorphisms.

We can give an explicit expression for this positive definite function on
$\Diff_c(X)$ since $F T(\psi) F $ can be written as a matrix relative
to the basis $\{f_n\}$ as $\left[  \ip{T(\psi) f_i}{f_j} \right]$.

The method of proof is the use the notion of generalized characters \cite{O1}
which are described fully in section \ref{sec:mainthm}.
This approach will give another proof of the positive definiteness of $\phi_{F,f}$.

We note that the representations $\pi_{F,f_+}$ or $\pi_{F,f_-}$
where $f_-(z) = 1-\beta + \beta z$, respectively
$f_+(z) = (1+ \alpha - \alpha z)^{-1}$,  are given in terms of the quasi-free
representations of the CAR, respectively, CCR algebras.  See \cite{B1,SV}.

\section{$N/V$ Limit and the Orthogonal Group}
\label{sec:nvlimit}

We include this section because it is hard to locate these results and they
are needed in our discussion in section \ref{sec:asymptotic}.

Let $L'$ be a real Hilbert space with a non-zero
vector $\omega,$ then the subgroup $K'$ of $O(L')$ of orthogonal
transformation $W$ such that $W \omega  = \omega$ can be
identified with the orthogonal group of $L' \ominus  \langle \omega \rangle$.
The irreducible spherical functions $\phi$ for $(O(L'), K')$ are
classified as 

\begin{equation} 
\phi^{(n)}(W) = 
\frac{1} {   \|   \omega   \|^{2n}     }
( S^n (W) \omega^{ \otimes n} , \omega^{\otimes n}),
\label{eq:sph}
\end{equation}

\noindent
where $S^n (W)$ denotes the $n$-th symmetric power of $W$,
which  is a classic result of I.J. Schoenberg.

More generally,  for an increasing sequence
$L_1 \subset L_2 \subset \dots$ of infinite dimensional
real Hilbert spaces whose union has completion $L$,
 we set $\Omega_k = f_1 + f_2 + \cdots + f_k$,
with $ \| \Omega_k \|^2 = V_k$.
Here, $ \{ f_j \} $ is an orthogonal system (not necessarily
orthonormal) such that $f_j \in L_j \ominus L_{j-1},$
$(L_0 = \{0\}$).
As usual, set $\displaystyle {\cal {O}} = \lim_\to O(L_k)$ and let
$\displaystyle {\cal K} = \lim_\to K_k,$
where $K_k$ consists of all elements in $O(L_k)$ that
fixed the vector $\Omega_k$.
Then the method of asymptotic limits given in Ol'shanskii \cite{O3}, theorem 22.10, gives the 
following result:

\begin{prop}
The irreducible spherical functions of the pair
$({\cal {O}},{ \cal K})$ have the form:
$\displaystyle
E_\lambda (W) = e^{ \lambda ([W - I] \Omega, \Omega)}, \,
\lambda \in \R^+ \cup \{0\},
$
where $\Omega$ is the formal vector $f_1+f_2 + \dots$.
$E_\lambda $ is given as the limit of irreducible
spherical functions of the pairs $(O(L_k), K_k),$
where
$\displaystyle E_\lambda = \lim_{k \to \infty} \, \phi^{(N_k)}_k$
if 
$\displaystyle \lambda = \lim_{k \to \infty} \, \frac{ N_k }{ \| \Omega_k \|}$
and where $\phi^{(N_k)}_k$ is the spherical function with index
$N_k$ of the group ${\cal {K}}_k$ given in equation \ref{eq:sph}.
\end{prop}

 If we restrict these irreducible spherical functions $E_\lambda$ to $\Diff_c(X)$
for representations constructed from an increasing sequence ${\cal F}$ of subsets
$X_1 \subset X_2 \subset \dots$ with $f_k = \chi_{X_k} - \chi_{X_{k-1}}$
(take $X_0 = \emptyset$),
so $\Omega_k = \chi_{X_k}$ and its norm $\| \Omega_k \|$ gives the volume $V_k$.
We find that

\[
\lim_{k \to \infty} \,
\phi^{(N_k)}_k ( T( \psi)) =  \exp[( [ T(\psi) - I] \Omega, \Omega)]
=
\exp \left[ \int_X \{ \det \psi(x) \}^{1/2} -1  \, dx \right],
\]

\noindent
which is the functional for the free Boson gas given in  \cite{G1, M,VGG}.

\section{Generalized Characters}

We recall that for a pair $(L,K)$,
where $L$ is a topological group with closed
subgroup $K$, that a $K$-central positive definite
function on $L$ is one that is invariant under
conjugation by $K$. The set of all such
normalized $K$-central positive definite
functions forms a convex set whose extreme
points  are called
\textsl {generalized characters} whose corresponding representations
are semifinite factor representations of type I or II.
Their importance for infinite dimensional
classical groups and their basic results were given  by Ol'shanskii  \cite{O1}.

\section{Complete Positivity and Positive Definite Functions}
\label{sec:cp}

Let $\chi$ be a holomorphic normalized positive definite function on the unitary
group $U(H)$, where $H$ is an infinite dimensional separable Hilbert space,
continuous in the strong operator topology.  Then $\chi$ is the restriction to $U(H)$
of the non-linear state of $B(H)$ (in the sense of Arveson \cite{A}) of the form

\[
\sum_{n=0}^\infty c_n \Tr [ A_n \pi_0^{\otimes n}( \cdot) ],
\]

\noindent
where
$A_n$ is a positive self-adjoint operator on $H^{\otimes n}$ with trace 1 and
$c_n \geq 0$ with $\sum_{n=0}^\infty c_n =1$, see \cite{P}.
 Further, such non-linear states are
completely positive, that is, if $x_{ij} \in B(H)$ with $[x_{ij}] \geq 0$, then
$[ \chi(x_{ij})] \geq 0$.  Further, following Ol'shanskii \cite{O2}, $\chi$ has a natural
extension, written
$\widehat{\chi}$,  to $\Gamma(H)$, the $*$-semigroup of contractions on $H$;
if $\pi_{\hat \chi}$ is the corresponding
$*$-representation of $\Gamma(H)$, then
$\{ \pi_{\hat \chi} (\Gamma(H)) \}'' = \{ \pi_{\chi}(U(H)) \}''$.

We let $\U = \lim_\to U(H_N)$,
where $H_0 \subset H_1 \subset H_2 \subset \cdots$ are all separable
infinite dimensional Hilbert spaces, as in \ref{sec:limitgroup}.  Again any holomorphic positive definite
function $\chi$ on ${\cal U}$ extends to $\Gamma = \lim_\to \Gamma(H_N)$
with
$\{ \pi_{\chi} ( \U ) \}'' = \{ \pi_{\hat \chi} ( \Gamma ) \}''$.

As in section \ref{sec:limitgroup}, we let $H$ have the orthonormal basis 
$\{e_j\} _{j=1}^\infty \cup \{ f_j\}_{j=1}^\infty$.
Set $H_0 = H_+ $ be the closed span of $\{e_j\} _{j=1}^\infty$
and
$H_-$ the closed span of $\{f_j\} _{j=1}^\infty$.
So, $H = H_+ \oplus H_-$. 
Set $H_N = H_0 \oplus \langle f_1, f_2, \dots, f_N \rangle$
and
$U(\infty) = \lim_\to U(  \langle f_1, f_2, \dots, f_N \rangle)$,
$U(2 \infty) = \lim_\to U( \langle e_1, \dots, e_N, f_1, \dots, f_N \rangle)$.
If $F$ denotes the projection of $H$ onto $H_-$ and $\Phi$ the compression map
 $\Phi(x) = FxF$, then $\Phi$ is a (linear) completely positive map.

Through $\Phi$, we obtain a canonical correspondence between positive definite
functions on $U(\infty)$ and $U(2 \infty)$ given as follows:

If $\chi$ is a holomorphic positive definite function on $U(\infty)$, let
$\widehat{\chi}$ denote its  extension  to $*$-semigroup $\Gamma(\infty)$ as a completely
positive function. Then $\widehat{\chi} \circ \Phi$ on $U(2\infty)$
is positive definite since it is the composition of completely positive
functions and a completely positive function is automatically positive definite.
Through continuity,
we can interpret $\widehat{\chi}\circ \Phi$ as a positive definite function on $\U$ as well.

If $\chi$ is, moreover,  a holomorphic finite
character on $U(\infty)$, we find that $\widehat{\chi} \circ \Phi$ are $U(\infty)$-conjugate invariant
positive definite functions on $U(2\infty)$.  In
section \ref{sec:mainthm} we find that they are, in fact,
generalized characters, that is, extremal in the convex set of $U(\infty)$-invariant 
positive definite functions.

\section
{The Asymptotic Character Formula as an $N/V$  Limit}
\label{sec:asymptotic}

As an application of the ideas of section \ref{sec:cp},  we  extend the asymptotic representation
of finite characters to generalized characters and relate it the $N/V$-limit.

The asymptotic character formula  \cite{B2, VK}
implies that for a fixed finite holomorphic character $\chi$ and
for a fixed $W \in U(\Co^{V_0})$ (here $V_0$ is an integer)
there exists a sequence of normalized characters
$\tilde{ \chi}_V$ of $U(\Co^V)$ such that 
$\tilde{\chi}_V (W) \to \chi(W)$.
In fact, such a sequence can be chosen relative to the
statistics of the signature $\lambda_V$.
Let $r_j ( \lambda_V)$ denote the length of its $j$-th row;
$c_j ( \lambda_V)$ of its $j$-th column; $ | \lambda_V| $ denote
the sum of the entries of the signature. Then
$\tilde { \chi}_V$ converges if and only if

\begin{equation}
\lim_{V \to \infty}
\frac{r_j ( \lambda_V)}{  V} = \alpha_j, \,
\lim_{V \to \infty}
\frac{ c_j ( \lambda_V)}{ V} = \beta_j, \,
\textrm{ and  } \,
\lim_{V \to \infty}
\frac{  | \lambda_V |}{ V} = \gamma.
\label{eq:limits}
\end{equation}

\noindent
where the parameters $\gamma, \alpha_j$, and $\beta_j$ determine
the mereomorphic function $f \in {\cal {M}}$ that gives
$\chi(W) = \det[ f(W)]$.

Since characters are central functions, we may view this
limiting process as taking place over the set of eigenvalues
$S^1 \times \cdots \times S^1$
instead of $U(\Co^{V_0})$.
Further, we know that $\tilde{\chi}_V$ is a sequence
of uniformly bounded holomorphic polynomials in $V_0$
complex variables
and $\chi$ itself is a meromorphic function in
$V_0$ complex variables.

Therefore, we also have pointwise convergence of
the holomorphic extensions of the normalized
characters; that is,
$\tilde{ \chi}_V$ to $\chi$ over the set
$D^1 \times \cdots \times D^1,$
where $D^1$ denotes the unit disk in the complex
plane. In particular, we find that
for fixed $W \in \Gamma(\Co^{V_0})$,
the set of contractions, that
$\tilde { \chi}_V (W) \to \det[ f(W)]$,
where $f = f_\chi$ is the meromorphic function of
one complex variable that defines the finite character $\chi$.

Finally, we let $\widehat \chi \circ \Phi$ be a generalized
character of $\U$.   Consider a fixed element $W \in \U$,
so that $W \in U(H_{V_0})$.
Then $\Phi(W) \in \Gamma(\Co^{V_0})$.
But we know already that there exists a sequence
of normalized characters $\tilde {\chi}_V$ such that
$\tilde { \chi}_V ( \Phi(W) ) \to \widehat \chi ( \Phi (W))$.

For easy comparison with the $N/V$ limit, we use the notation in
section \ref{sec:nvlimit} so let
$\U = \lim_{\to} U( L^2(X_V))$.
Further we consider the special case
$\chi(W)= \det[f_+(W)]$,
for $f_+(z) = {1}/{(1+\beta_1 - \beta_1 z)}$.
Then the  corresponding finite factor representation is an inductive
limit of symmetric algebras and can be interpreted in terms of
a quasi-free representation of the Weyl (CCR) algebra, see \cite{B1}.
Let $\beta_1 = \lim_{ V \to \infty} \frac{N_V}{V}$,
where $V$ and $N_V$ are positive integers which correspond to the degree of the 
symmetric tensor 
power and the dimension of the representation space, respectively.
(Note that the other limits in equation (\ref{eq:limits})  are all zero.)
Set 

\[ \tilde {\chi}_V (W) 
= \Tr[ S^{N_V} ( \Phi_V (W))] / \textrm{dim}(S^{N_V})(  \Co^{V}  )),
\]

\noindent
where $W \in U(L^2(X_V))$ and
 $V$ gives the volume of $X_V$. 
Then $\tilde {\chi}_V (W) \to \widehat \chi(\phi(W))$;  in particular,
$\tilde{\chi}_V ( T( \psi)) \to \widehat \chi ( \phi ( T( \psi)))$.

Hence, the translation of $V$ in the asymptotic formula to the $N/V$ limit
is to change its role as the dimension of
$\Co^{V}$ to the volume of the $X_V$ for the Hilbert space $L^2(X_V)$.

\section{Proof of the Main Theorem \ref{thm:main} }
\label{sec:mainthm}

The method of proof is to identify the functions $\phi_{F,f}$ as extreme
points in the convex set of all generalized characters. Then by \cite{O1}
we know that the corresponding representation is semifinite and factor.
Although these functions were introduced in \cite{SV}, there are no general
results concerning the structure of their corresponding representations.

\begin{lemma} \label{lemma:l1}
Let $ L = \lim_\to U( 2 \infty) 
  = \lim_\to U(\langle  e_1, \dots, e_N, f_1, \dots, f_N  \rangle),$
with subgroups
$K_1 = U(\langle  e_1, e_2, \dots\rangle) \times I$,
$K_2 = I \times U(\langle f_1, f_2, \dots  \rangle)$.
Then every holomorphic generalized character $\phi$ of the
pair $(L,K_2)$ such that $\phi | K_1 =1$ can be
written uniquely as $\phi = \widehat \chi \circ \Phi$,
where $\widehat \chi$ is the holomorphic extension of a 
holomorphic finite character $\chi$ from $U(\infty)$
to $\Gamma(\infty)$.
\end{lemma}

\pf
Let $\phi$ be a generalized character of $(L,K_2)$ whose restriction
to $K_1$ is trivial.
Set $\phi_{2n} = \phi | U(2n)$.
Then we make the claim that

\[
\phi_{2n} ( \cdot )
= \sum_\alpha 
c^{(2n)}_\alpha \Tr[ \pi_\alpha^{(2n)}(E_n) \pi_\alpha^{(2n)} 
( \cdot )],
\]

\noindent
where $E_n$ is the projection of $\langle e_1, \dots , e_n, f_1, \dots, f_n\rangle$
onto $\langle f_1, \dots, f_n \rangle$.

To see this, write $\phi_{2n}$ as 
$\sum_\alpha \Tr[P_\alpha^{(2n)} \pi_\alpha^{(2n)}( \cdot)]$,
where $P_\alpha^{(2n)}$ is a positive operator on the representation space 
$H( \pi_\alpha^{(2n)} )$ and 
$\sum_\alpha \Tr [ P_\alpha^{(2n)} ]=1$,
by section 1.
The unitary branching law:
$
\pi_\alpha^{(2n)} |  U(n) \times  U(n) \simeq
   \sum_\beta  \pi_\beta^{(n)} \times \pi_{\alpha / \beta}^{(n)}
$
together with the condition $ \phi_{2n} | U(n) \times I$
is trivial implies that the support of the operator
$P_\alpha^{(2n)}$ must be contained in the projection
onto the subspace where 
$   \pi_\beta^{(n)} \times \pi_{\alpha / \beta}^{(n)} $
acts as  $ \pi_0^{(n)} \times \pi_\alpha^{(n)}.$
In particular, $P_\alpha^{(2n)} = c_\alpha^{(2n)} \pi_\alpha^{(2n)}(E_n),$
for some $c_\alpha^{(2n)} \geq 0.$
Further, the signatures $\alpha$ that contribute to $\phi_{2n}$
have at most $n$ non-zero entries.

As usual, we let $\pi_\alpha^{(2n)}$ denote the holomorphic extension of
$\pi_\alpha^{(2n)}$ from $U( 2n)$ to $\Gamma( 2n)$.
Then, by the centrality of the trace, we have:

\begin{equation}
\phi_{2n} ( W )
= \sum_\alpha c_\alpha^{(2n)} \Tr[ \widehat \pi_\alpha^{(2n)}(E_n W)]
=
\sum_\alpha c_\alpha^{(2n)} \Tr[ \widehat \pi_\alpha^{(n)}(\Phi( W))].
\label{eq:above}
\end{equation}

\noindent
By (\ref{eq:above}), we find that $\phi_{2n}$ is uniquely determined by its restriction
to $I \times U(n)$.  Moreover, the map $\phi \mapsto \phi | I \times U(\infty)$
preserves convex combinations.

Conversely, let $\widehat \chi$ denote the extension to $\Gamma(\infty)$ 
of a holomorphic finite character $\chi$ of $U(\infty).$ 
Then $\widehat \chi \circ \Phi$ is trivial on $K_1$ and conjugate invariant
by elements of $K_2$. Further,  this correspondence is the
inverse of the map given above: $\phi \mapsto \phi | K_2$,
since $\widehat \chi \circ \Phi | K_2 = \chi$.
\smallskip

\begin{lemma} \label{lemma:l2}
Given an orthonormal basis $\{ e_j\}$ for $H_{+}$  
such that, if 
 $S_N = \{ j : e_j \in H_N \},$ 
$\{ e_j : j \in S_N \} \cup \{f_j\}_{j=1}^N $ is an orthonormal basis
for $H_N,$
where $\{ f_N\}$ is the orthonormal basis  for $H_-$ as above,
we form
 $\displaystyle U(2 \infty) = \lim_\to U( \langle   e_1,\dots, e_N, f_1, \dots, f_N  \rangle)$.
Then
$U( 2 \infty)$ is dense in ${\cal {U}}$.
\end{lemma}

\pf \
Let $W \in \U$ be given so $W \in U(H_{N_0})$,
for some $N_0.$
Let $\epsilon > 0$ be given, together with
vectors $\xi_1, \dots, \xi_k$.
Then we need to find $W' \in U(\langle  e_1, \dots, e_{n'}, f_1, \dots, f_{n'}  \rangle)$,
for some $n'$
such that $ \|  W \xi_j - W'\xi_j \| < \epsilon$,
for $1 \leq j \leq k$.

We consider for each $j$ the orthonormal expansion
$
\xi_j = \sum_{n=1}^\infty c_n^{(j)} e_n + \sum_{n=1}^{N_0} d_n^{(j)} f_n.
$
Choose the index $n'$ so that 
$ \|   \sum_{n=n'}^\infty c_n^{(j)} e_n  \| < \epsilon/2$,
for $1 \leq j \leq k.$
We note that since $\xi_j \in H_{N_0}$
the coefficients $c_n^{(j)} = 0,$ if $e_n \notin H_{N_0}$.
For this reason, we find that
$W ( \sum_{n=1}^{n'} c_n^{(j)} e_n) \in H_{N_0}$.
Hence, there exists $W' \in U(\langle   e_1,\dots,e_{n'},  f_1,\dots,f_{n'}  \rangle)$
such that
$\|   W \xi_j - W' \xi_j'  \|  < \epsilon/2,$ for $1 \leq j \leq n$,
where $\xi_j' = \sum_{n=1}^{n'} c_n^{(j)} e_n$.
This is the desired $W'$,
and the Lemma is proven.
\smallskip

\begin{prop} \label{prop:main}
(1)
The  convex sets $C_1$ and $C_2$ are affinely isomorphic,
where $C_1$ is the set of all holomorphic positive definite
functions on $\U$ that are $I \times U(\infty)$-conjugate invariant
and whose restriction to $(U(H_+) \cap \U) \times I$ is trivial,
and where 
$C_2$ is the set of all central holomorphic positive definite
functions on $U(\infty)$.

\noindent
(2)
Let $\phi$ be any holomorphic generalized character from
$C_2$.
Then there exists a unique holomorphic finite character $\chi$
of $U(\infty)$ such that $\phi = \widehat \chi \circ \phi$.
\end{prop}

\pf \
Let $\phi$ be a generalized character of the pair
$(\U, I \times U(\infty)).$
Then $\phi | U(2 \infty)$ is a generalized character
as well of $(U(2 \infty), I \times U(\infty))$,
for a version of $U(2 \infty)$ that satisfies 
Lemma \ref{lemma:l2}.
Further, since $U(2\infty)$ is dense in $\U$,
$\phi$ is uniquely determined by its restriction
to $U(2 \infty)$.
By Lemma \ref{lemma:l1}, $\phi = \widehat \chi \circ \Phi$
on $U(2 \infty)$.
But $\widehat {\chi} \circ \Phi$ is continuous on
$U(2 \infty)$ relative to the $\U$-topology.
Hence, $\phi = \widehat {\chi} \circ \Phi$ on
$\U$ uniquely.
\smallskip

\begin{cor} \label{cor:above}
Let $\chi$ denote a holomorphic finite character of $U(\infty)$.
Then $\widehat {\chi} \circ \phi$ is a generalized character of $\U$,
which generates either an irreducible representation or
a type $\textrm{II}_\infty$ factor representation.
In particular, the generalized character is irreducible
if and only if its restriction to $U(\infty)$ is a power
of the determinant; it is type  $\textrm{II}_\infty$
if and only if its restriction to $U(\infty)$ is
type $\textrm{II}_1$.
\end{cor}

\pf \
We make use of the identification of the generalized
characters as spherical functions of the pair
$(L \times K, K \times K)$   \cite{O1}.
If $f$ is the generalized character with corresponding
cyclic representation $T = \pi \times \pi'$ of
$L \times K,$
then $T$ is irreducible,  $\pi'$ is finite factor,
and $\pi$ and $\pi'$ generate each other commutants.
Hence, $\pi$ is
irreducible $\iff \pi'$ is irreducible.
This occurs if and only if $\pi'$ is equivalent
to a power of the determinant.
Otherwise, $\pi'$ is always equivalent to
a
$\textrm{type II}_1$ factor representation.
So, $\pi$ must be type II as well.
Since $U(H_N),$ so $\U$ itself, has no infinite dimensional 
finite factor representations, we find that
$\pi$ is type $\textrm{II}_\infty$ if and only if  $\pi'$ is
type $\textrm{II}_1$.

\medskip

We now turn to the proof of the main theorem \ref{thm:main}.
By the discussion in Section \ref{sec:limitgroup}, the von Neumann algebras
generated by $\pi_{F,f}( \U)$ and $\pi_{F,f}( \Diff_c(X))$ are identical.
Hence, the main theorem follows from \ref{prop:main} and its corollary.

\medskip

We close with indicating two possible further extensions of this work.

\begin{enumerate}

\item It would be interesting to incorporate other classes of representations into the
scheme presented in this paper. For example, the tensor product of the semifinite
representations of \ref{thm:main} with the standard representations of the
diffeomorphism group that come from restricting the tame representations of
$U(L^2(X))$ as well as the tensor product of the spherical representations discussed in
\cite{O3}. This would be a much larger family than the representations presented in
\cite{VGG}.  

\item Further,  we would like to explore  the connection of our construction
of the generalized characters by the condition expectation operator with the Rieffel
approach to induced representations.  See \cite{Gr,L}.

\end{enumerate}

\medskip

\noindent
\textsc{Department of Mathematics and Computer Science,
Drexel University,
Philadelphia, PA 19104}

\noindent
\textit{E-mail address:} \texttt{rboyer@mcs.drexel.edu}

\medskip

\noindent
\textsc{keywords}:{Diffeomorphism group, factor representation, unitary group} 

\noindent
\textsc{Mathematics Subject Classification}:{Primary 22E65, 81R10; Secondary 46L99}

\end{document}